\documentclass[11pt]{article}
\newtheorem{theorem}{Theorem}[section]
\newtheorem{lemma}[theorem]{Lemma}
\newtheorem{corollary}[theorem]{Corollary}
\newtheorem{proposition}[theorem]{Proposition}
\newtheorem{df}[theorem]{Definition}
\newenvironment{definition}{\begin{df} \rm}{\end{df}}

 \usepackage{amsmath,amssymb,amscd}
\newcommand{\pf}[1]{\trivlist \item[\hskip\labelsep\it #1\ ]}
\newcommand{\varpf}[1]{\trivlist \item[\hskip\labelsep\sc #1:]}
\newcommand{\qedbox}{$\rlap{$\sqcap$}\sqcup$}
\newcommand{\qed}{\qquad \qedbox \endtrivlist}
\newcommand{\varqed}{\hfill \rule{0.6em}{0.6em} \endtrivlist}
\newenvironment{proof}{\pf{Proof}}{\qed}

\newenvironment{remark}{\pf{Remark}}{\endtrivlist}
\newenvironment{remarks}{\pf{Remarks} 
   \begin{enumerate}}{\end{enumerate} \endtrivlist}

\newenvironment{examples}{\pf{Examples} 
   \begin{enumerate}}{\end{enumerate} \endtrivlist}
\newenvironment{items}{
  \begin{enumerate} 
                    
  }{\end{enumerate}}

\newenvironment{classification}{\noindent\small 2000 {\it Mathematics Subject
Classification\/}:}{\vskip 12pt}

\newcommand{\defiff}{\quad:\Longleftrightarrow\quad}

\newcommand{\posints}{{\mathbb N}}

\newcommand{\tensor}{\otimes}

\newcommand{\ptensor}{\hat{\otimes}}

\newcommand{\wtensor}{\check{\otimes}}

\newcommand{\cstar}{{C^\ast}}

\newcommand{\im}{{\operatorname{Im}}}
\newcommand{\id}{{\mathrm{id}}}

\newcommand{\rad}{{\operatorname{rad}}}

\newcommand{\A}{{\mathfrak A}}

\begin{document}
\title{Banach space properties \\ forcing a reflexive, amenable Banach algebra to be trivial}
\author{\it Volker Runde}
\date{}
\maketitle
\begin{abstract}
It is an open problem whether an infinite-dimensional amenable Banach algebra exists whose underlying
Banach space is reflexive. We give sufficient conditions for a reflexive, amenable Banach algebra to be
finite-dimensional (and thus a finite direct sum of full matrix algebras). If $\A$ is a reflexive, amenable Banach algebra such 
that for each maximal left ideal $L$ of $\A$ (i) the quotient $\A / L$ has the approximation property and (ii) the canonical map
from $\A \wtensor L^\perp$ to $(\A / L) \wtensor L^\perp$ is open, then $\A$ is finite-dimensional. As an application,
we show that, if $\A$ is an amenable Banach algebra whose underlying Banach space is an ${\cal L}^p$-space with
$p \in (1,\infty)$ such that for each maximal left ideal $L$ the quotient $\A / L$ has the approximation property, then
$\A$ is finite-dimensional.
\end{abstract}
\section*{Introduction}
The question of whether a reflexive, amenable Banach algebra $\A$ has to be trivial, i.e.\ if it is necessarily of the form
\[
  \A \cong {\mathbb M}_{n_1} \oplus \cdots \oplus {\mathbb M}_{n_k}
\]
with $n_1, \ldots, n_k \in \posints$, was apparently raised for the first time by J.\ E.\ Gal\'e, T.\ J.\ Ransford, and
M.\ C.\ White in \cite{GRW}. In the same paper, they showed that every reflexive, amenable Banach algebra all of whose 
primitive ideals have finite codimension must be trivial. Generalizing this result, B.\ E.\ Johnson proved in \cite{Joh3} 
that every reflexive, amenable whose maximal left ideals are complemented has to be trivial. This implies, for instance,
that the underlying Banach space of an amenable Banach algebra cannot be a Hilbert space (\cite{GLW}). For a recent refinement
of Johnson's result, see \cite{Zha}. 
\par
In general, complemented subspaces are rather the exception than the rule in functional analysis (except, of course, in
the Hilbert space setting), which makes theorems involving hypotheses of complementedness often less applicable 
than one would like them to have. In this paper, we establish criteria for the triviality of reflexive, amenable Banach
algebras which avoid complemented subspaces in their hypotheses. Instead, we require certain Banach spaces properties
such as the (bounded) approximation property.
\par
In the $\cstar$-algebra situation, the connection between amenability and Banach space properties is well established:
An amenable $\cstar$-algebra is nuclear (\cite{Con,BP}), and thus has the (metric) approximation property (\cite{CE}).
In \cite{LLW}, this fact was used to prove that $\mbox{$\ell^\infty$-}\bigoplus_{n=1}^\infty {\mathbb M}_n$ cannot
be amenable. In \cite{Run1}, we showed that any reflexive, amenable Banach algebra with the approximation property
such that every bounded linear map from $\A$ into $\A^\ast$ is compact, has to be trivial (in fact, it is sufficient
that the almost periodic functionals separate the points of $\A$; see \cite{Run2}). In conjunction with Pitt's theorem
(\cite[p.\ 208]{Koe}), this implies that the underlying Banach space of an amenable Banach algebra cannot be $\ell^p$ 
with $p > 2$. 
\par
In this paper, we shall first establish a general theorem --- involving no complemented subspaces in its hypotheses ---
which gives two sufficient conditions for a reflexive, amenable Banach algebra to be trivial. We shall then apply this
theorem to Banach algebras whose underlying Banach spaces are ${\cal L}^p$-spaces in the sense of \cite{LR} with
$p \in (1,\infty)$.
\section{Preliminaries}
Amenable Banach algebras were introduced by B.\ E.\ Johnson in \cite{Joh1} in terms of a cohomological triviality
condition. In this paper, however, we require a characterization of amenable Banach algebras given in \cite{Joh2}.
\par
We write $\tensor$ for the algebraic tensor product of Banach spaces, $\ptensor$ for the projective tensor product,
and $\wtensor$ for the injective tensor product; our standard source on Banach space tensor products is
\cite{DF}. If $\A$ is a Banach algebra, then $\A \ptensor \A$ (and $\A \wtensor \A$,
as well) is a Banach $\A$-bimodule with the operations
\[
  a \cdot (x \tensor y) := ax \tensor y \quad\text{and}\quad (x \tensor y) \cdot a := x \tensor ya
  \qquad (a, x, y \in \A).
\]
We use $\Delta \!: \A \ptensor \A \to \A$ to denote the linear operator induced by the algebra multiplication, i.e.\
$\Delta(a \tensor b) = ab$ for $a,b \in \A$.
\par
An {\it approximate diagonal\/} for a Banach algebra $\A$ is a bounded net $( m_\alpha )_\alpha$ in
$\A \ptensor \A$ such that
\begin{equation} \label{diag1}
  a \cdot m_\alpha - m_\alpha \cdot a \to 0 \qquad (a \in \A)
\end{equation}
and
\begin{equation} \label{diag2}
  a \Delta m_\alpha \to a \qquad (a \in \A).
\end{equation}
As was shown in \cite{Joh2}, a Banach algebra $\A$ is amenable if and only if it has an approximate diagonal.
If an amenable Banach algebra has an approximate diagonal bounded by $\lambda \geq 1$, we say that $\A$
is $\lambda$-amenable.
\par
For any two Banach spaces $E$ and $F$, we write ${\cal F}(E,F)$ for the bounded finite rank operators from
$E$ into $F$; by ${\cal A}(E,F)$, we denote the approximable operators from $E$ into $F$, i.e.\ the norm 
closure of ${\cal F}(E,F)$. We simply write ${\cal F}(E)$ and ${\cal A}(E)$ insteady of ${\cal F}(E,E)$ and 
${\cal A}(E,E)$, respectively. A Banach space $E$ is said to have the {\it approximation property\/}, if there is
a net $( T_\alpha)_\alpha$ in ${\cal F}(E)$ such that $T_\alpha \to \id_E$ uniformly on compact subsets of $E$.
If this net $( T_\alpha)_\alpha$ can be chosen to be bounded in norm by some $\lambda \geq 1$, we say that
$( T_\alpha)_\alpha$ has the {\it $\lambda$-approximation property\/}; the $1$-approximation property is also
called {\it metric approximation property\/}. If a Banach space has the $\lambda$-approximation property for
some $\lambda \geq 1$, it is said to have the {\it bounded approximation property\/}. All these approximation
properties can also be defined in terms of nets in ${\cal A}(E)$.
\par
In general, the approximation property is not inherited by closed subspaces or quotients (see \cite{LT}). 
For certain amenable Banach algebras (and the appropriate sub- and quotient structures), however, such hereditary 
properties exist. We say that a Banach algebra $\A$ has compact multiplication if, for each $a \in \A$, the
operators $L_a$ and $R_a$ of multiplication by $a$ from the left and from the right, respectively, are compact.
\begin{proposition} \label{comp}
Let $\A$ be a semisimple $\lambda$-amenable Banach algebra with compact multiplication which has the approximation property.
Then each closed, one-sided ideal of $\A$ and each quotient of $\A$ by such an ideal has the $\lambda$-approximation
property. In particular, $\A$ itself has the $\lambda$-approximation property.
\end{proposition}
\begin{proof}
Let $( m_\alpha )_\alpha$ be an approximate diagonal for $\A$ bounded by $\lambda$. Then, by (\ref{diag1}) and 
(\ref{diag2}), 
$( \Delta m_\alpha )_\alpha$ is a bounded approximate identity for $\A$ bounded by $\lambda$. From the structure result
\cite[Corollary 3.3]{LRRW}, it follows that $\A$ has a bounded approximate identity $( e_\beta )_{b \in {\mathbb B}}$
bounded by $\lambda$ such that, for each $\beta \in \mathbb B$, the operators $L_{e_\beta}$ and $R_{e_\beta}$
belong to ${\cal F}(\A)$.
\par
Let $L$ be a closed left ideal of $\A$ (right ideals can be treated analoguously). Then 
$(L_{e_\beta} |_L)_{\beta \in {\mathbb B}}$ is a net in
${\cal F}(L)$ bounded by $\lambda$ such that $L_{e_\beta} |_L \to \id_L$ pointwise and hence uniformly
on compact subsets of $L$. A similar argument shows that $\A / L$ has the $\lambda$-approximation property.
\end{proof}
\begin{remark}
If we only require $\A$ to have a left approximate identity bounded by $\lambda \geq 1$ and 
that $\A$ has compact left multiplication, we can still show that, for each closed left ideal $L$ of
$\A$, both $L$ and $\A / L$ have what one might call the compact $\lambda$-approximation property 
(see \cite[I, p.\ 94]{LT} for the definition of the compact approximation property).
\end{remark}
\section{A general theorem}
In this section, we shall prove a general theorem which establishes that a reflexive, amenable Banach
algebra $\A$ has to be trivial provided that two additional demands are satisfied. The first of
these requirements is that every quotient of $\A$ modulo a maximal left ideal has to have the approximation
property, and the second one is that a certain tensor product of maps between injective tensor products
has to be open (or, equivalently, surjective).
\par
The technical heart of our argument is the following lemma (if $F$ is a subspace of a Banach space $E$, we write
$F^\perp$ for its annihilator in $E^\ast$):
\begin{lemma} \label{main}
Let $\A$ be an amenable Banach algebra, and let $L \neq \A$ be a closed left ideal of $\A$ such that:
\begin{items}
\item $\A / L$ has the bounded approximation property;
\item if $\pi_L \!: \A \to \A / L$ is the canonical surjection and, $\pi_L \tensor \id_{L^\perp} \!:
\A \wtensor L^\perp \to (\A / L) \wtensor L^\perp$ is open.
\end{items}
Then there is $e \in \A^{\ast\ast} \setminus L^{\ast\ast}$ such that $L e = \{ 0 \}$.
\end{lemma}
\begin{proof}
Let $( m_\alpha )_{\alpha \in {\mathbb A}}$ be an approximate diagonal for $\A$. For each $\alpha \in \mathbb A$,
$m_\alpha$ has the form
\[
  m_\alpha = \sum_{n=1}^\infty a_n^{(\alpha)} \tensor b_n^{(\alpha)},
\]
where we may suppose that $\sum_{n=1}^\infty \| a_n^{(\alpha)} \| < \infty$ and $\lim_{n\to \infty} \| b_n^{(\alpha)} \| = 0$.
It follows that, for each $\alpha \in \mathbb A$, the set $\{ \pi_L(b_n^{(\alpha)}) : n \in \posints \}$ is relatively 
(norm) compact. Since $\A / L$ was supposed to have the bounded approximation property, there is $\lambda \geq 1$ such
that, for each $\alpha \in \mathbb A$ and for each $\epsilon > 0$, there is $S_{\alpha, \epsilon} \in
{\cal F}(\A/L)$ with $\| S_{\alpha,\epsilon} \| \leq \lambda$ and
\[
  \| S_{\alpha, \epsilon} (\pi_L(b^{(\alpha)}_n)) - \pi_L(b^{(\alpha)}_n) \| 
  < \frac{\epsilon}{1+ \sum_{n=1}^\infty \| a_n^{(\alpha)} \|}
  \qquad (n \in \posints).
\]
Since $( \A / L)^\ast \cong L^\perp$, we may identify ${\cal A}(\A/L)$ with $(\A / L) \wtensor L^\perp$
(\cite[4.2]{DF}). By (ii) there is a constant $C \geq 1$ such that, for each ${\bf x} \in (\A / L) \wtensor L^\perp$,
there is ${\bf y} \in \A \wtensor L^\perp$ with $\| {\bf y} \| \leq C \| {\bf x} \|$ and 
$(\pi_L \tensor \id_{L^\perp})({\bf y}) = {\bf x}$.
We may thus find, for each $\alpha \in \mathbb A$ and $\epsilon > 0$, an operator $T_{\alpha,\epsilon}
\in {\cal A}(\A/L,\A) \cong \A \wtensor L^\perp$ with $\| T_{\alpha,\epsilon} \| \leq \lambda C$ such that
$\pi_L \circ T_{\alpha,\epsilon} = S_{\alpha,\epsilon}$. Define
\[
  e_{\alpha,\epsilon} := \sum_{n=1}^\infty a_n^{(\alpha)} T_{\alpha,\epsilon}(\pi_L(b_n^{(\alpha)}))
  \qquad (\alpha \in {\mathbb A}, \, \epsilon > 0).
\]
It is easy to see that $( e_{\alpha, \epsilon} )_{\alpha,\epsilon}$ is a bounded net in $\A$ and thus has
a $w^\ast$-accumulation point $e \in \A^{\ast\ast}$. Without loss of generality suppose that
$e = \mbox{$w^\ast$-}\lim_{\alpha,\epsilon} e_{\alpha,\epsilon}$. We show that 
$\pi_L^{\ast\ast}(e) \neq 0$. We have:
\begin{eqnarray*}
  \pi_L^{\ast\ast}(e) & = & \mbox{$w^\ast$-}\lim_{\alpha,\epsilon} \pi_L 
  \left( \sum_{n=1}^\infty a_n^{(\alpha)} T_{\alpha,\epsilon}(\pi_L(b_n^{(\alpha)})) \right) \\
  & = & \mbox{$w^\ast$-}\lim_{\alpha,\epsilon}
  \sum_{n=1}^\infty a_n^{(\alpha)} \cdot \pi_L(T_{\alpha,\epsilon}(\pi_L(b_n^{(\alpha)}))) \\
  & = & \mbox{$w^\ast$-}\lim_{\alpha,\epsilon} 
  \sum_{n=1}^\infty a_n^{(\alpha)} \cdot S_{\alpha,\epsilon}(\pi_L(b_n^{(\alpha)})), 
  \qquad\text{by the choice of $T_{\alpha,\epsilon}$}, \\
  & = & \mbox{$w^\ast$-}\lim_\alpha \sum_{n=1}^\infty a_n^{(\alpha)} \cdot \pi_L(b_n^{(\alpha)}),
  \qquad\text{by the choice of $S_{\alpha,\epsilon}$}, \\
  & = & \mbox{$w^\ast$-}\lim_\alpha \pi_L \left( \sum_{n=1}^\infty a_n^{(\alpha)} b_n^{(\alpha)} \right) \\
  & \neq & 0, \qquad\text{by (\ref{diag2})}.
\end{eqnarray*}
\par
Let $a \in L$. Then
\begin{equation} \label{zero}
  \sum_{n=1}^\infty a_n^{(\alpha)} T_{\alpha,\epsilon}(\pi_L(b_n^{(\alpha)}a)) = 0
  \qquad (\alpha \in {\mathbb A}, \, \epsilon > 0)
\end{equation}
since $L$ is a left ideal. We then have:
\begin{eqnarray*}
  ae & = & 
  \mbox{$w^\ast$-}\lim_{\alpha,\epsilon} \sum_{n=1}^\infty a a_n^{(\alpha)} T_{\alpha,\epsilon}(\pi_L(b_n^{(\alpha)})) \\
  & = & \mbox{$w^\ast$-}\lim_{\alpha,\epsilon} \sum_{n=1}^\infty a_n^{(\alpha)} T_{\alpha,\epsilon}(\pi_L(b_n^{(\alpha)}a)),
  \qquad\text{by (\ref{diag1})}, \\
  & = & 0, \qquad\text{by (\ref{zero})}.
\end{eqnarray*}
This completes the proof.
\end{proof}
\begin{remark}
We do not know if condition (ii) of Lemma \ref{main} is automatically satisfied: There are examples of Banach spaces
$E$ and $X$ and of a closed subspace $F$ of $E$ such that $\pi_F \tensor \id_X \!: E \wtensor X \to (E/F) \wtensor X$
is not open (\cite[4.3]{DF}). These examples, however, do not involve any Banach algebra structure, let alone 
amenability.
\end{remark}
\begin{theorem} \label{refcor}
Let $\A$ be a reflexive, amenable Banach algebra such that for each maximal left ideal $L$ of $\A$:
\begin{items}
\item $\A / L$ has the approximation property;
\item if $\pi_L \!: \A \to \A / L$ is the canonical surjection, $\pi_L \tensor \id_{L^\perp} \!:
\A \wtensor L^\perp \to (\A / L) \wtensor L^\perp$ is open.
\end{items}
Then $\A$ is trivial.
\end{theorem}
\begin{proof}
Let $L$ be a maximal left ideal of $\A$.
Since $\A$ is reflexive, so is $\A / L$. Hence, by \cite[16.4, Corollary 4]{DF}, $\A / L$ has the metric 
approximation property. By Lemma \ref{main}, there is thus an element $e \in \A \setminus L$ such that $Le = \{ 0 \}$. 
Since $L$ was arbitrary, it follows that $\A / \rad(\A)$ is a modular right annihilator algebra in the sense
of \cite[8.4.6 Definition]{Pal}. By \cite[8.4.14 Proposition]{Pal}, this means that $\A / \rad(\A)$ is
finite-dimensional. It follows from standard results on the splitting of exact sequences (compare the proof
of \cite[Proposition 2.3]{Run1}) that $\A$ must be trivial.
\end{proof}
\par
Let $\A$ be a Banach algebra, and let $L$ be a {\it complemented\/} left ideal in $\A$. Then $\pi_L \tensor \id_{L^\perp} \!:
\A \wtensor L^\perp \to (\A / L) \wtensor L^\perp$ has a bounded right inverse, and thus, in particular, is open. 
As a special case of Theorem \ref{refcor} we therefore obtain \cite[Proposition 1.10]{GLW}:
\begin{corollary}
Let $\A$ be an amenable Banach algebra whose underlying Banach space is a Hilbert space. Then $\A$ is trivial.
\end{corollary}
\section{${\cal L}^p$-spaces for $p \in (1,\infty)$}
Of course, it would be rather unsatisfactory if all what Theorem \ref{refcor} could accomplish were to recover
the Hilbert space case. In this section, we consider Banach algebras whose underlying Banach space is
an ${\cal L}^p$-space for $p \in (1,\infty)$. For $p \in [1,\infty]$, the ${\cal L}^p$-spaces were introduced
in \cite{LP}, and were later investigated in greater detail in \cite{LR}. All $L^p$-spaces are ${\cal L}^p$-spaces.
For $p \in (1,\infty)$, the ${\cal L}^p$-spaces are reflexive
(\cite[Theorem I(i)]{LR}), and are isomorphic to complemented subspaces of
$L^p$-spaces (\cite[Theorem I(ii)]{LR}).
\par
We shall work, however, in a more abstract framework:
\begin{definition} \label{PQdef}
Let $E$ be a Banach space, and let $\lambda \geq 1$. Then: 
\begin{items}
\item $E$ has property $(P_\lambda)$ if, for each finite-dimensional subspace $F$ of $E$, there is a finite-dimensional
subspace $X$ of $E$ containing $F$ such that there is a projection of norm not greater than $\lambda$ onto $X$.
\item $E$ has property $(Q_\lambda)$ if, for each closed subspace $F$ of $E$ with finite codimension, there is
a closed subspace $X$ of $F$ with finite codimension in $E$ such that the canonical surjection $\pi_X \!: E \to E / X$
has a bounded right inverse of norm not greater than $\lambda$. 
\end{items}
\end{definition}
\begin{remark}
Trivially, every Banach space with property $(P_\lambda)$ for some $\lambda \geq 1$ has the $\lambda$-approximation 
property.
\end{remark}
\begin{examples}
\item Every Hilbert space has both property $(P_1)$ and $(Q_1)$.
\item If $E$ is an ${\cal L}^p$-space with $p \in [1,\infty]$, then $E$ has property 
$(P_\lambda)$ for some $\lambda \geq 1$ (\cite[Theorem  III(c)]{LR}).
\end{examples}
\par
There is a natural relation between properties $(P_\lambda)$ and $(Q_\lambda)$:
\begin{lemma} \label{dual}
Let $E$ be a Banach space, let $\lambda \geq 1$, and suppose that $E^\ast$ has property $(P_\lambda)$. Then
$E$ has property $(Q_{\lambda(1+\epsilon)})$ for every $\epsilon > 0$.
\end{lemma}
\begin{proof}
Let $F$ be a closed subspace of $E$ with finite codimension. Then $F^\perp$ is a finite-dimensional subspace of
$E^\ast$. Since $E^\ast$ has property $(P_\lambda)$, there is a finite-dimensional subspace $Y$ of $E^\ast$ containing
$F^\perp$ such that there is a projection $P \!: E^\ast \to Y$ with $\| P \| \leq \lambda$. Let $X = {^\perp Y}$.
Then $X$ is a closed subspace of $E$ with finite codimension such that $E / X \cong E^{\ast\ast} / X^{\ast\ast} \cong
Y^\ast$ isometrically. Let $\pi_X \!: E \to E/X$ be the canonical projection. With the appropriate identifications
in place, it follows that $P^{\ast\ast} \!: E/X \to E^{\ast\ast}$ is a right inverse of $\pi_X^{\ast\ast}$ with
$\| P^{\ast\ast} \| \leq \lambda$.
\par
Let $\epsilon > 0$, and let $\phi_1, \ldots, \phi_n \in E^\ast$ be a Hamel basis of $Y$. Since $\im\, P^{\ast\ast}$
is finite-dimensional, it follows from the principle of local reflexivity (see, e.g., \cite[6.6]{DF}) that there is a 
linear map $T \!: \im \, P^{\ast\ast} \to E$ with $\|T \| \leq 1 +\epsilon$ such that $T_{E \cap \im\, P^{\ast\ast}} =
\id_{E \cap \im\, P^{\ast\ast}}$ and
\begin{equation} \label{locref}
  \langle T P^{\ast\ast}x, \phi_j  \rangle = \langle \phi_j, P^{\ast\ast} \rangle \qquad
  (x \in E^{\ast\ast}, \, j =1, \ldots, n).
\end{equation}
Let $Q := T P^{\ast\ast}$. Then (\ref{locref}) implies that $Q \!: E/ X \to E$ is a right inverse of $\pi_X$, and
certainly $\| Q \| \leq \lambda(1+\epsilon)$.
\end{proof}
\begin{remark}
If $E$ is reflexive, the application of the local reflexivity principle in the proof of Lemma \ref{dual} becomes
unnecessary, so that property $(P_\lambda)$ for $E^\ast$ implies property $(Q_\lambda)$ for $E$.
\end{remark}
\par
The relevance of property $(Q_\lambda)$ in connection with condition (ii) of Lemma \ref{main} becomes apparent
in our next lemma:
\begin{lemma} \label{l2}
Let $E$ be a reflexive Banach space which has property $(Q_\lambda)$ for some $\lambda\geq 1$, let $F$ be 
a closed subspace of $E$, and let $\pi_F \!: E \to E/F$ be the canonical map. Then, for any Banach space $X$, 
$\pi_F \tensor \id_X \!: E \wtensor X \to (E/F) \wtensor X$ is open. 
\end{lemma}
\begin{proof}
We consider $\pi_F \tensor \id_X \!: E \tensor X \to (E/F) \tensor X$. We shall prove that, for each ${\bf y} 
\in (E/F) \tensor X$ there is ${\bf x} \in E \tensor X$ with $(\pi_F \tensor \id_X)({\bf x}) = {\bf y}$ and
$\| {\bf x} \| \leq \lambda \| {\bf y} \|$, where the norm is the injective norm. Certainly, ${\bf y} 
\in (E/F) \tensor Z$ for some finite-dimensional subspace  $Z$ of $X$. Since the injective norm respects subspaces
isometrically (\cite[4.3, Proposition]{DF}), we may suppose without loss of generality that $\dim X < \infty$.
Note that under this hypothesis, $E \tensor X = E \wtensor X$ is a reflexive Banach space.
\par
Let ${\cal Y}$ be the collection of all closed subspaces of $E$ containing $F$ and having finite codimension. Defining
\[
  Y_1 \prec Y_2 \defiff Y_1 \supset Y_2 \qquad (Y_1, Y_2 \in {\cal Y})
\]
Since $E$ has property $(Q_\lambda)$, there is, for each $Y \in {\cal Y}$, a closed subspace $\tilde{Y}$ of
$Y$ with finite codimension in $E$ such that the canonical surjection $\pi_{\tilde{Y}} \!: E \to E / \tilde{Y}$ has
a right inverse of norm not greater than $\lambda$.
\par
Let ${\bf y} \in (E/F) \tensor X$. Then, by the foregoing, there is, for each $Y \in \cal Y$, an element
${\bf x}_Y \in E \tensor X$ with $\| {\bf x}_Y \| \leq \lambda \| {\bf x} \|$ such that 
$(\pi_Y \tensor \id_X)({\bf x}_Y) = (\tilde{\pi}_Y \tensor \id_X)({\bf y})$ (here, 
$\tilde{\pi}_Y \!: E/F \to E/Y$ is the map induced by the
canonical surjection $\pi_Y \!: E \to E/Y$). The net $( {\bf x}_Y )_{Y \in {\cal Y}}$ is bounded in 
$E \tensor X$ and thus has a weak accumulation point ${\bf x}$. It follows that
\[
  (\pi_Y \tensor \id_X)({\bf x}) = (\tilde{\pi}_Y \tensor \id_X)({\bf y}) \qquad (Y \in {\cal Y}),
\]
and hence that $(\pi_F \tensor \id_X)({\bf x}) = {\bf y}$.
\end{proof}
\par
As an immediate consequence of Theorem \ref{refcor} and Lemma \ref{l2}, we obtain:
\begin{theorem} \label{thm2}
Let $\A$ be a reflexive, amenable Banach algebra such that:
\begin{items}
\item $\A / L$ has the approximation property for each maximal left ideal $L$ of $\A$;
\item $\A$ has property $(Q_\lambda)$ for some $\lambda \geq 1$.
\end{items}
Then $\A$ is trivial.
\end{theorem}
\par
For ${\cal L}^p$-spaces, we obtain:
\begin{corollary} \label{Lp}
Let $\A$ be an amenable Banach algebra whose underlying Banach space is an ${\cal L}^p$-space with $p \in (1,\infty)$.
Then one of the following holds:
\begin{items}
\item $\A$ is trivial;
\item there is a maximal left ideal $L$ of $\A$ such $\A / L$ lacks the approximation property.
\end{items}
\end{corollary}
\begin{proof}
By \cite[Theorem I(i)]{LR}, $\A$ is reflexive, and by \cite[Theorem III(a)]{LR}, $\A^\ast$ is an ${\cal L}^q$-space,
where $q^{-1} + p^{-1} = 1$. As we already noted after Definition \ref{PQdef}, $\A^\ast$ has property
$(P_\lambda)$ for some $\lambda \geq 1$. Hence, by the remark following Lemma \ref{l2}, $\A$ has property 
$(Q_{\lambda})$ as well. Consequently, if $\A / L$ has the approximation property for each maximal left ideal
$L$ of $\A$, Theorem \ref{thm2} forces $\A$ to be trivial.
\end{proof}
\begin{remarks}
\item Let $p \in (1,\infty) \setminus \{ 2 \}$, and let $E$ be an ${\cal L}^p$-space. Then it follows from 
\cite[Theorem I(iv)]{LR} and either \cite[I, Theorem 2.d.6]{LT} (for the case $p < 2$) or
\cite[II, Theorem 1.g.4]{LT} (for the case $p > 2$) that $E$ has a quotient which lacks the approximation property. 
For $p \neq 2$, we thus cannot {\it a priori\/} rule out case (ii) in Corollary \ref{Lp}.
\item Both \cite[I, Theorem 2.d.6]{LT} and \cite[II, Theorem 1.g.4]{LT} are pure Banach space results, i.e.\
they involve no Banach {\it algebra\/} structure whatsoever. Hence, it may well be that case (ii) of
Corollary \ref{Lp} does never occur. In view of Proposition \ref{comp}, we raise the following formal question: 
If $\A$ is an amenable Banach algebra which has the approximation property, does $\A / L$ have the approximation
property for each closed left ideal $L$ of $\A$?
\end{remarks}
\vfill
\begin{tabbing}
{\it Address\/}: \= Department of Mathematical Sciences \\
\> University of Alberta \\
\> Edmonton, Alberta \\
\> Canada T6G 2G1 \\[\medskipamount]
{\it E-mail\/}: \> {\tt runde@math.ualberta.ca}
\end{tabbing} 
\end{document}